\documentclass[11pt]{article}

  \RequirePackage[dvips]{hyperref}

\usepackage{amsmath,amssymb,amscd,amsthm,mathrsfs}
\usepackage{graphicx}
\usepackage{graphics}
\usepackage{float}

\usepackage{amsmath,amsfonts}
\usepackage{mathrsfs}
\usepackage{graphicx}
\newtheorem{theorem}{\bf Theorem}[section]

\newtheorem{lemma}{\bf Lemma}[section]

\def \Liminf{\mathop{\underline{\lim}}\limits}

\def\1{\mbox{1\hspace{-.25em}I}}
\def\Ex{{\bf E}}
\def\Pb{{\bf P}}
\let\bb\mathbb
\def\UU{{\bb U}}
\def\BB{{\bb B}}
\def\CC{{\bb C}}
\def\NN{{\bb N}}

\def\sgn{{\rm sgn}}

\begin{document}

\title{{\bf On Parameter Estimation of Threshold
Autoregressive Models}}

\author{Ngai Hang Chan and Yury A. Kutoyants\\
\small Chinese University of Hong Kong and Universit\'e du Maine}
\date{}

\maketitle

\begin{abstract}
This paper studies the threshold estimation of a TAR model
when the underlying threshold parameter is a random variable. It is
shown that the Bayesian estimator is consistent and its limit distribution is
expressed in terms of a limit likelihood ratio. Furthermore, convergence of moments of
the estimators is also established.   The limit
distribution can be computed via explicit simulations from which testing
and inference for the threshold parameter can be conducted.  The obtained
results are illustrated with numerical simulations. \\


\end{abstract}

\section[short title]{Introduction}
Since the publication of the seminal treatise of Tong \cite{To90},
the field of nonlinear time series has been receiving considerable attention
in the literature.  Today, nonlinear time series has been widely applied to
subjects such as ecology, engineering, chaos, finance and econometrics.  From a
statistical perspective, nonlinear time series also furnishes an exciting
platform for nonstandard statistical inference both parametrically
and nonparametrically.  For a comprehensive survey on some of these recent
developments, see Fan and Yao \cite{FY03}.

Among many different developments in nonlinear time series, estimation and
testing of the threshold parameter constitute one of most challenging tasks.
One of the main reasons of the difficulty arises from the fact that tricky and
nonstandard asymptotic techniques are required to handle the threshold
estimation, see Chan \cite{Ch93}, Hansen \cite{Han97} and \cite{Han00}.
A comprehensive theory for this type of problems seems to be lacking from
the literature so far, however.

On the other hand, a relative complete theory for the statistical inference
for diffusion processes in continuous time is available, see for example
Kutoyants \cite{Kut98} and \cite{Kut04}.  In particular, these two books
demonstrate that both the maximum likelihood and the Bayesian approaches to
diffusion processes can be put under a general context and an asymptotic theory
can be developed, albeit to its non standard nature.

One of the main purposes of this paper is to make use of this general theory
and apply it to the nonlinear time series context.  Related contributions
to continuous time ARMA and threshold ARMA models can be found, for example, in
Brockwell \cite{Brock94}, Chan and Tong \cite{CT87}, Stramer, Brockwell and
Tweedie \cite{SBT96}, and Tong \cite{To90} and the references therein.

Although likelihood inference for the threshold parameter of nonlinear time series
was considered by Chan \cite{Ch93} and Hansen \cite{Han97} previously, the asymptotic
machineries  employed were of special nature which cannot be easily generalized to other
situations.  For further background on likelihood tests of non-linearity,
see Li and Li \cite{LL08}.  From a Bayesian perspective, Geweke and Teuri \cite{GT93}
considered a Bayesian threshold AR model and derived the posterior distribution of
the threshold parameter.  However, a detailed description of the asymptotic properties of
the Bayesian estimator and its moment convergence were lacking.

By incorporating the developments in diffusion, this paper illustrates a general
methodology to tackle both the maximum likelihood and Bayesian estimation problems from
which simulations can be efficiently conducted.  Moreover, the proposed
approach is sufficiently transparent and can be easily adopted to other
nonlinear time series context.

A second but equally important goal of this study is to develop an implementable
scheme for simulating and computing the limit likelihood statistics.  By linking
the integral equation of the underlying invariant density of the nonlinear
time series and the intensity of the limiting compound Commission process, one can
compute the form of the limiting likelihood explicitly.  To the best of our
knowledge, this has never been conducted before and results obtained in this paper
can greatly enhance the inference for the threshold parameter of a nonlinear time
series and extend its applications.

This paper is organized as follows.  Background introduction together with the
statement of the problem and the main result are given in Section 2.  Section 3
consists of simulations. Section 4 discusses the extension to cover the
usual one-sided threshold setting while conclusions and possible extensions are
given in Section 5.

\section[short title]{Main result}

Consider the model
\begin{equation}
\label{2-1}
X_{j+1}=\rho_1 \,X_j\,\1_{\left\{\left|X_j\right|<\vartheta
\right\}}+\rho_2 \,X_j\,\1_{\left\{\left|X_j\right|\geq \vartheta
\right\}}+\varepsilon _{j+1}, \qquad j=0,\ldots, n-1,
\end{equation}
where $\varepsilon _j$ are i.i.d. ${\cal N}\left(0,\sigma^2 \right)$,  $\rho
_1\not =\rho _2$  and
$\left|\rho_2\right| <1$.  Note that model~(\ref{2-1}) appears to be
different from the standard setting, where the threshold is usually
partitioned as $\1_{\left\{X_j<\vartheta \right\}}$ and $\1_{\left\{X_j \geq \vartheta
\right\}}$.  We choose the current setting because it is more general
and mathematically more convenient.  Our results can be easily
extended to encompass the standard setting as demonstrated in Section 4.1.  We suppose that $\sigma^2>0,\rho
_1,\rho _2 $ are known and $\vartheta \in \Theta =\left(\alpha ,\beta
\right)$ is the unknown threshold parameter. Our goal is to estimate $\vartheta  $ from
observations $X^n=\left(X_0,X_1,\ldots,X_n\right)$ and  to
describe the asymptotic behavior of the estimators as $n\rightarrow \infty $.
Recall that $\left(X_j\right)_{j\geq 1}$ is geometrically mixing (see Chen
and Tsay \cite{CT91}) and denote its stationary density function
by $f\left(\vartheta ,x\right)$, see also Fan and Yao \cite{FY03}.

In this paper, we consider both the maximum likelihood and Bayesian approaches.
Recall that the likelihood function is written as
\begin{align*}
L\left(\vartheta ,X^n\right)&=f_0\left(X_0\right)
\left(\frac{1}{\sqrt{2\pi \sigma ^2}}\right)^n \exp\Bigl\{-\frac{1}{2\sigma
^2} \sum_{j=0}^{n-1}\left(X_{j+1}-\rho_1
\,X_j\,\1_{\left\{\left|X_j\right|\leq \vartheta \right\}}\right.\\
&\qquad\qquad\qquad\qquad \left. -\rho_2
\,X_j\,\1_{\left\{\left|X_j\right|>\vartheta \right\}}\right)^2\Bigr\}\,,
\end{align*}
and the maximum likelihood estimator (MLE) $\hat \vartheta _n$ is defined by the equation
\begin{equation}
\label{mle}
\sup_{\vartheta \in\Theta }L\left(\vartheta ,X^n\right)=\max
\left[L\left(\hat\vartheta_n+ ,X^n\right),L\left(\hat\vartheta_n-
,X^n\right)\right] .
\end{equation}
If this equation has many solutions, we can, for example, call the MLE to be the
value which is at the center of the gravity. Note that
the function $L\left(\vartheta ,X^n\right),\vartheta \in \Theta $ has jumps at
the points
$$
\vartheta _l=\left|X_j\right|\in \Theta ,\qquad l=1,\ldots, L,
$$
where $L\leq n$.  Clearly, if $\Theta =R$, then $L=n$.

To apply the Bayesian approach, suppose that the unknown parameter is a random
variable with a known prior density $p\left(\theta \right), \theta \in \Theta$,
which is continuous and positive.  Using
the quadratic loss function, the Bayesian estimator (which
minimizes the mean squares error) is the conditional mathematical expectation
\begin{equation}
\label{be}
\tilde \vartheta _n=\int_{\alpha }^{\beta }\theta \,p\left(\theta \right)L\left(\vartheta
,X^n\right)\,{\rm
d}\theta
=\frac{\int_{\alpha }^{\beta }\theta \,p\left(\theta  \right)L\left(\vartheta
,X^n\right)\,{\rm d}\theta}{\int_{\alpha }^{\beta }p\left(\theta  \right)L\left(\vartheta
,X^n\right)\,{\rm d}\theta}.
\end{equation}

Properties of the least squares estimator (LSE) of $\vartheta$ were studied in Chan \cite{Ch93}.
The LSE coincides with the MLE for Gaussian $\varepsilon _j$.  We therefore recall
properties of MLE and compare them with properties of the Bayesian estimators.

First, introduce the stochastic process
\begin{align*}
Z\left(u\right)= \left\{\begin{array}{ll}
             &\exp\left\{  -\frac{\rho ^2\vartheta ^2}{2\sigma ^2}\;\;
N_+\left(u\right)-\frac{\rho\, \vartheta}{\sigma ^2}
\,\sum_{l=0}^{N_+\left(\;u\,\,\,\,\right)} \varepsilon_l^+
\right\},\qquad   u\geq 0 , \\
             &\exp\left\{ -\frac{\rho ^2\vartheta ^2}{2\sigma ^2}\;\;
N_-\left(-u\right)-\frac{\rho\, \vartheta}{\sigma ^2}
\,\sum_{l=0}^{N_-\left(-u\right)} \varepsilon_l^-  \right\},\quad u\leq 0,
             \end{array}
             \right.
\end{align*}
where $N_+\left(\cdot \right)$ and $N_-\left(\cdot \right)$ are two independent
Poisson processes of intensities $\lambda _+=\lambda _-=2f\left(\vartheta,\vartheta
\right)$ ($f(\vartheta, x)$ is the stationary density function of $X_j$) and
$\varepsilon_r^+, \varepsilon_l^- $ are independent Gaussian
${\cal N}\left(0,\sigma ^2\right)$ random variables. It is easy to see that
$$
Y_+\left(u\right)=\rho^2\,
\vartheta^2\,N_+\left(u\right)+2\rho\, \vartheta
\,\sum_{l=0}^{N_+\left(u\right)} \varepsilon_l^+ =\sum_{l=0}^{N_+\left(u\right)}\left[\rho^2\,
\vartheta^2+2\rho\, \vartheta \,\varepsilon_l^+ \right],\qquad u\geq 0
$$
$$
Y_-\left(u\right)=\rho^2\,
\vartheta^2\,N_-\left(u\right)+2\rho\, \vartheta
\,\sum_{l=0}^{N_-\left(u\right)} \varepsilon_l^- =\sum_{l=0}^{N_-\left(u\right)}\left[\rho^2\,
\vartheta^2+2\rho\, \vartheta \,\varepsilon_l^- \right],\qquad u\geq 0
$$
are compound Poisson processes.

The random process $Z\left(\cdot \right)$ is piecewise constant and as a result,
the points $ u^*$ of the maximum of the process $Z\left(\cdot \right)$ is defined
by
$$
\sup_{u}Z\left(u\right)=Z\left( u^*\right)\,,
$$
where
$$
\hat u_m< u^*<\hat u_M.
$$
Here $\hat{u}_m$ and $\hat{u}_M$ are two consecutive events of the process
$N_+\left(\cdot \right)$, or of the process $N_-\left(\cdot \right)$, or they are
respectively the first event of $N_+\left(\cdot \right)$ and $N_-\left(\cdot \right)$.
Simulated realizations of $Z\left(\cdot \right)$ are given in Section 3.
The center of gravity of the interval is given by the point
\begin{equation}
\label{choice}
\hat u=\frac{u_m+u_M}{2}.
\end{equation}
Such a choice of $\hat u$ is explained in Section 4 below.
It follows from the result of
Chan \cite{Ch93} that the MLE $\hat \vartheta _n$
is consistent and
$$
n\left(\hat \vartheta _n-\vartheta \right)\Longrightarrow \hat u.
$$
\noindent
Introduce the random variable
$$
\tilde u=\frac{\int_{}^{}u\,Z\left(u\right)\,{\rm
d}u}{\int_{}^{}Z\left(u\right)\,{\rm d}u}.
$$
\noindent
The main result is the following theorem.

\noindent
\begin{theorem}
\label{T1}
The Bayesian estimator $\tilde \vartheta _n$ constructed by the observations
$X^n$ of the threshold autoregressive process is consistent, the normalized
difference $n\left(\tilde \vartheta _n-\vartheta \right) $ converges in
distribution :
\begin{equation}
\label{2-2}
n\left(\tilde \vartheta _n-\vartheta \right)\Longrightarrow \tilde u
\end{equation}
and for any $p>0$
\begin{equation}
\label{2-3}
\lim_{n\rightarrow \infty }\Ex_\vartheta \left|n\left(\tilde \vartheta
_n-\vartheta \right)\right|^p =\Ex_\vartheta \left|\tilde u\right|^p.
\end{equation}
\end{theorem}

\noindent
{\bf Proof.} The proof of this theorem is based on the general result by
Ibragimov and Khasminskii \cite{IH81}, Theorem 1.10.2. To apply it we study the
normalized likelihood ratio process
\begin{align*}
Z_n\left(u\right)=\frac{L\left(\vartheta+\frac{u}{n}
,X^n\right)}{L\left(\vartheta ,X^n\right)} ,\qquad u\in
\UU_n=\left[n\left(\alpha -\vartheta \right),n\left(\beta -\vartheta
\right)\right],
\end{align*}
where $\vartheta $ is the true value.
Recall the main steps. Write the Bayesian estimator
$\left(\theta_u =\vartheta +\frac{u}{n}\right)$ as
\begin{align*}
\tilde \vartheta_n&=\frac{\int_{\alpha }^{\beta }\theta \,p\left(\theta  \right)L\left(\theta
,X^n\right)\,{\rm d}\theta}{\int_{\alpha }^{\beta }p\left(\theta  \right)L\left(\theta
,X^n\right)\,{\rm d}\theta}=\vartheta +\frac{1}{n}\frac{\int_{\UU_n }^{}u
\,p\left(\theta_u  \right)L\left(\theta_u
,X^n\right)\,{\rm d}u}{\int_{\UU_n }p\left(\theta_u  \right)L\left(\theta_u
,X^n\right)\,{\rm d}u}\\
&=\vartheta +\frac{1}{n}\frac{\int_{\UU_n }^{}u
\,p\left(\theta_u  \right)\frac{L\left(\theta_u
,X^n\right)}{L\left(\vartheta
,X^n\right) }\,{\rm d}u}{\int_{\UU_n }p\left(\theta_u  \right)\frac{L\left(\theta_u
,X^n\right)}{L\left(\vartheta
,X^n\right)}\,{\rm d}u}=\vartheta +\frac{1}{n}\frac{\int_{\UU_n }^{}u
\,p\left(\theta_u  \right)Z_n\left(u\right)\,{\rm d}u}{\int_{\UU_n }
p\left(\theta_u  \right)Z_n\left(u\right)\,{\rm d}u}\,.
\end{align*}

Suppose that we proved the convergence of the process
$Z_n\left(\cdot \right)$ to the process $Z\left(\cdot \right)$ providing the
convergence of these integrals. Then
$$
n\left(\tilde \vartheta_n-\vartheta \right)\Longrightarrow \frac{\int_{}^{}u \,
Z\left(u\right)\,{\rm d}u}{\int_{}
Z\left(u\right)\,{\rm d}u}=\tilde u.
$$
 This convergence together with an estimate on the large deviations of the
 tails of the process $Z_n\left(\cdot \right)$ allow us to prove the convergence
 of the moments \eqref{2-3}.

Now check the conditions of the Theorem 1.10.2 in \cite{IH81}. We
need to prove
\begin{enumerate}
\item  the convergence of the finite dimensional distributions of
$Z_n\left(\cdot \right)$ to the finite dimensional distributions of
$Z\left(\cdot \right)$, that is,
\begin{equation} \label{fdd}
Z_n\left(\cdot \right) \rightarrow Z\left(\cdot \right) \ \ \mbox{f.d.d.},
\end{equation}
\item to establish the  estimate:
\begin{equation}
\label{2-4}
\Ex_\vartheta
\left[Z_n^{1/2}\left(u_2\right)-Z_n^{1/2}\left(u_1\right)\right]^2\leq
C\,\left|u_2-u_1\right|\,,
\end{equation}
\item and to establish the estimate: for any $M>0$
\begin{equation}
\label{2-5}
\Ex_\vartheta
Z_n^{1/2}\left(u\right)\leq \frac{C_M}{\left|u\right|^M}\,.
\end{equation}
\end{enumerate}

The convergence of finite-dimensional distributions follows from the Proposition 2 of
\cite{Ch93}.  Instead of repeating a technical argument as in \cite{Ch93}, we offer a
different intuitive (but rigorous) explanation as follows.  Rewrite the process \eqref{2-1} as
\begin{equation}
\label{2-6}
X_{j+1}=\rho_1 \,X_j+\rho \,X_j\,\1_{\left\{\left|X_j\right|\geq \vartheta
\right\}}+\varepsilon _{j+1}, \qquad j=0,\ldots, n-1,
\end{equation}
where we use $\1_{\left\{\left|X_j\right|< \vartheta
\right\}}=\1_{}-\1_{\left\{\left|X_j\right|\geq \vartheta
\right\}}$ and denote $\rho =\rho _2-\rho_1$.

Put $ Z_n\left(u\right)=\exp\left\{-\frac{1}{2\sigma
^2}Y_n\left(u\right) \right\}$ and study the process $Y_n\left(u\right)$ for
positive values of $u$.
\begin{align*}
Y_n\left(u\right)& =\sum_{j=0}^{n-1}\left[\left(X_{j+1}-\rho_1 \,X_j-\rho
\,X_j\,\1_{\left\{\left|X_j\right|>\vartheta+\frac{u}{n}
\right\}}\right)^2\right.\\
&\quad\qquad \qquad \left. -\left(X_{j+1}-\rho_1 \,X_j-
\rho\,X_j\,\1_{\left\{\left|X_j\right|>\vartheta \right\}}\right)^2\right]\\
&=\sum_{j=0}^{n-1}\left(
\rho\,X_j\,\left[\1_{\left\{\left|X_j\right|>\vartheta\right\}}-
\1_{\left\{\left|X_j\right|>\vartheta+\frac{u}{n}\right\}}\right] \right)\\
&\qquad \times \left(2X_{j+1}-2\rho_1 \,X_j-
\rho\,X_j\,\left[\1_{\left\{\left|X_j\right|>\vartheta\right\}}+
\1_{\left\{\left|X_j\right|>\vartheta+\frac{u}{n}\right\}}\right]  \right).
\end{align*}
Note that
$$
\1_{\left\{\left|X_j\right|>\vartheta\right\}}-
\1_{\left\{\left|X_j\right|>\vartheta+\frac{u}{n}\right\}}=\1_{\left\{\vartheta
<\left|X_j\right|\leq \vartheta+\frac{u}{n}\right\}}.
$$
Hence,
\begin{align}
Y_n\left(u\right)&=\sum_{j=0}^{n-1}\rho\,X_j \left[ 2  X_{j+1}-2\rho_1 \,X_j-\rho\,X_j
\right]\;\1_{\left\{\vartheta
<\left|X_j\right|\leq\vartheta+\frac{u}{n}\right\}}\nonumber\\
&=\sum_{j=0}^{n-1}\left(\rho ^2\,X_j^2+2\rho \,X_j\,\varepsilon _{j+1} \right)\;\1_{\left\{\vartheta
<\left|X_j\right|\leq\vartheta+\frac{u}{n}\right\}}  .
\label{lr}
\end{align}

Next introduce another process
\begin{equation}
\label{lro}
Y_n^\circ\left(u\right)=\sum_{j=0}^{n-1}\left(\rho ^2\vartheta ^2+ 2\rho
\,\vartheta\,\sgn\left(X_j\right)\,
\varepsilon _{j+1}\right) \;\1_{\left\{\vartheta
<\left|X_j\right|\leq\vartheta+\frac{u}{n}\right\}}
\end{equation}
and put $\1_{\left\{\vartheta
<\left|X_j\right|\leq\vartheta+\frac{u}{n}\right\}}=\1_{\left\{\BB_j\left(u\right)\right\}}$.
We show that this process is asymptotically equivalent to the process
$Y_n\left(u\right)$. We have
\begin{align*}
\Ex_\vartheta \left|Y_n\left(u\right)-Y_n^\circ\left(u\right) \right|&\leq \rho
^2\sum_{j=0}^{n-1} \Ex_\vartheta \left|X_j^2-\vartheta
^2\right|\1_{\left\{\BB_j\left(u\right)\right\}}\\
& +2\rho \sum_{j=0}^{n-1}
\Ex_\vartheta
\left|X_j-\vartheta \,\sgn\left(X_j\right)
\right| \left|\varepsilon _{j+1}\right|\1_{\left\{\BB_j\left(u\right)\right\}}.
\end{align*}
For the first term we write
\begin{align*}
\Ex_\vartheta \left|X_j^2-\vartheta
^2\right|\1_{\left\{\BB_j\left(u\right)\right\}}&=\int_{\vartheta \leq
\left|x\right|\leq \vartheta +\frac{u}{n}}^{}\left|x^2-\vartheta
^2\right|\,f\left(\vartheta ,x\right)\,{\rm d}x \\
&\leq 2\left(\vartheta +\frac{u}{n}\right)\left(\frac{u}{n}\right)^2\max_{\vartheta
<\left|x\right|<\vartheta+\frac{u}{n}}f\left(\vartheta ,x\right)\leq
C\,\left(\frac{u}{n}\right)^2.
\end{align*}
  The second term is (recall that $X_j$ and $\varepsilon _{j+1}$ are independent)
\begin{align*}
&\Ex_\vartheta
\left|X_j-\vartheta \,\sgn\left(X_j\right)
\right| \left|\varepsilon
_{j+1}\right|\1_{\left\{\BB_j\left(u\right)\right\}}\\
&\qquad\qquad\qquad =  \Ex_\vartheta
\left|X_j-\vartheta \,\sgn\left(X_j\right)
\right| \,\1_{\left\{\BB_j\left(u\right)\right\}}  \Ex_\vartheta \left|\varepsilon
_{j+1}\right|\\
&\qquad \qquad\qquad=\sqrt{\frac{2}{\pi
}}\;\sigma \; \Ex_\vartheta
\left|X_j-\vartheta \,\sgn\left(X_j\right)
\right| \,\1_{\left\{\BB_j\left(u\right)\right\}}\\
&\qquad \qquad\qquad= 2\sqrt{\frac{2}{\pi
}}\;\sigma \; \int_{\vartheta \leq
x\leq \vartheta +\frac{u}{n}}^{}\left|x-\vartheta
\right|\,f\left(\vartheta ,x\right)\,{\rm d}x \leq C\,\left(\frac{u}{n}\right)^2.\\
\end{align*}
Hence,
$$
\Ex_\vartheta \left|Y_n\left(u\right)-Y_n^\circ\left(u\right) \right|\leq
C\;\frac{u^2}{n}\longrightarrow 0
$$
for any fixed $u$.  Therefore, it is sufficient to study the limit distribution of the random
function $Y_n^\circ\left(u\right)$ and to show the convergence
\begin{equation}
\label{ll}
Y_n^\circ\left(u\right)\Longrightarrow
Y_+\left(u\right)=\sum_{l=0}^{N_+\left(u\right)}\left[\rho^2\,
\vartheta^2+2\rho\, \vartheta \,\varepsilon_l^+ \right].
\end{equation}

To see that the limit of $Y_n^\circ\left(u\right)$ is a
compound Poisson process, first note that the characteristic function
\begin{align}
\label{2-7}
\Phi\left(v\right)&=\Ex_\vartheta e^{{\rm i}vY_+\left(u\right)}=\Ex_\vartheta e^{{\rm
i}v\sum_{l=0}^{N_+\left(u\right)}\left[\rho^2\,
\vartheta^2+2\rho\, \vartheta \,\varepsilon_l^+ \right]  }\nonumber\\
&=\left.\Ex_\vartheta \Ex_\vartheta \left( e^{{\rm
i}v\sum_{l=0}^{N_+\left(u\right)}\left[\rho^2\,
\vartheta^2+2\rho\, \vartheta \,\varepsilon_l^+ \right]  }\right|{\cal
F}_{N_+}\right)\nonumber\\
&=\Ex_\vartheta  e^{\left[{\rm
i}v\rho^2\,
\vartheta^2-2v^2\rho^2\, \vartheta^2 \,\sigma ^2\right]N_+\left(u\right)
}\nonumber\\
&=\exp\left\{u\left(e^{ {\rm i}v\rho^2\,
\vartheta^2-2v^2\rho^2\, \vartheta^2 \,\sigma ^2 }-1\right)2f\left(\vartheta
,\vartheta \right)\right\},
\end{align}
where we denote ${\cal F}_{N_+} $ to be the $\sigma $-algebra related to the
Poisson  process and make use of the independence of $\varepsilon _l^+$ and
$N_+\left(\cdot \right)$.
The desired convergence will be proved if  the convergence of the
characteristic function of the process $Y_n^\circ\left(\cdot \right)$ to
\eqref{2-7} is established.

Fix $u>0$, then as $n\rightarrow \infty $ the band $\left[\vartheta
,\vartheta +\frac{u}{n}\right]$ becomes narrower and the events,
when
$\left|X_{j_l}\right|\in \left[\vartheta ,\vartheta +\frac{u}{n}\right]$, become
more rare. This means that the distance between two consecutive events
$\left|X_{j_l}\right|\in \left[\vartheta ,\vartheta +\frac{u}{n}\right]$ and
$\left|X_{j_{l+1}}\right|\in \left[\vartheta ,\vartheta +\frac{u}{n}\right]$ tends
to infinity. As the process $\left(X_j\right)_{j\geq 1}$ is geometrically
mixing, these events become asymptotically independent.  Under such circumstances,
the characteristic function
$\Phi_n\left(v\right)=\Ex_\vartheta e^{{\rm i}vY_n^\circ\left(u\right)}$ can be
calculated explicitly as
\begin{align*}
\Phi_n\left(v\right)&=\left.\Ex_\vartheta\left(\Ex_\vartheta
\exp\left\{\sum_{j=0}^{n-1}{\rm i}v\left(\rho ^2\vartheta ^2+ 2\rho
\,\vartheta\,\sgn\left(X_j\right)\,
\varepsilon _{j+1}\right)
\;\1_{\left\{\BB_j\left(u\right)\right\}}\right\}\right|{\cal F}_X\right)\\
&=\Ex_\vartheta\exp\left\{\sum_{j=0}^{n-1}\left({\rm i}v\rho ^2\vartheta ^2- 2v^2\rho^2
\,\vartheta^2\,\sigma ^2\right)
\;\1_{\left\{\BB_j\left(u\right)\right\}}\right\}\\
&=\left(\Ex_\vartheta\exp\left\{ \left({\rm i}v\rho ^2\vartheta ^2- 2v^2\rho^2
\,\vartheta^2\,\sigma ^2 \right) \;\1_{\left\{\BB_1\left(u\right)\right\}}
\right\}\right)^n .
\end{align*}
Further,
\begin{align*}
&\Ex_\vartheta\exp\left\{ \left({\rm i}v\rho ^2\vartheta ^2- 2v^2\rho^2
\,\vartheta^2\,\sigma ^2 \right) \;\1_{\left\{\BB_1\left(u\right)\right\}}
\right\}\\
&\quad=\int_{ }^{ }e^{\left({\rm i}v\rho ^2\vartheta ^2- 2v^2\rho^2
\,\vartheta^2\,\sigma ^2\right)\1_{\left\{\BB_1\left(u\right)\right\}}}f\left(\vartheta ,x\right){\rm d}x \\
&\quad =\left(\int_{-\infty }^{-\vartheta -\frac{u}{n} }+\int_{\vartheta +\frac{u}{n}
}^{\infty}+\int_{-\vartheta
}^{\vartheta }\right)f\left(\vartheta ,x\right){\rm d}x\\
&\qquad+ \left(\int_{-\vartheta -\frac{u}{n}}^{-\vartheta} +\int_{\vartheta
}^{\vartheta+\frac{u}{n}}\right) e^{{\rm i}v\rho ^2\vartheta ^2- 2v^2\rho^2
\,\vartheta^2\,\sigma ^2  }f\left(\vartheta ,x\right){\rm d}x\\
&\quad= 1-2\frac{u}{n}f\left(\vartheta ,\vartheta \right)+2\frac{u}{n}e^{{\rm i}v\rho ^2\vartheta ^2- 2v^2\rho^2
\,\vartheta^2\,\sigma ^2 }f\left(\vartheta ,\vartheta \right)+o\left(\frac{u}{n}\right).
\end{align*}
Hence,
\begin{align*}
\ln \Phi_n\left(v\right)&=n\ln\left( 1+ \frac{u}{n}\left( e^{{\rm i}v\rho
^2\vartheta ^2- 2v^2\rho^2 \,\vartheta^2\,\sigma ^2 }-1\right)2f\left(\vartheta
,\vartheta \right)+o\left(\frac{u}{n}\right)\right) \\
&\longrightarrow u\left( e^{{\rm i}v\rho
^2\vartheta ^2- 2v^2\rho^2 \,\vartheta^2\,\sigma ^2 }-1\right)2f\left(\vartheta
,\vartheta \right)=\ln\Phi\left(v\right)\,.
\end{align*}
That is, it coincides with \eqref{2-7} and as a result, (\ref{fdd}) is proved. \hfill $\Box$

To establish conditions \eqref{2-4} and \eqref{2-5} we need the following two
lemmas.
\noindent
\begin{lemma}
\label{L1} There exists a constant $C>0$ such that for all values
  $u_1,u_2\in\left(n\left(\alpha -\vartheta\right),n\left(\beta  -\vartheta
 \right) \right)$ we have the inequality
$$
\Ex_\vartheta\left(Z_n^{1/2}\left(u_2\right)-Z_n^{1/2}\left(u_1\right)\right)^2\leq
C\,\left|u_2-u_1\right|.
$$
\end{lemma}
\noindent
{\bf Proof.} The first step is
\begin{align*}
\Ex_\vartheta
\left(Z_n^{1/2}\left(u_2\right)-Z_n^{1/2}\left(u_1\right)\right)^2&=2-2\;\Ex_\vartheta
\left[Z_n\left(u_2\right)Z_n\left(u_1\right)
\right]^{1/2}\\
&=2-2\;\Ex_{\vartheta+\frac{u_1}{n}}
\left[\frac{Z_n\left(u_2\right)}{Z_n\left(u_1\right)} \right]^{1/2},
\end{align*}
where the measure is changed from $\Pb_\vartheta$ to
$\Pb_{\vartheta+\frac{u_1}{n}}$.
We have  ($u_2\geq u_1>0$)
\begin{align*}
&\Ex_{\vartheta+\frac{u_1}{n}}
\left[\frac{Z_n\left(u_2\right)}{Z_n\left(u_1\right)}
\right]^{1/2}\\
&\quad =
\Ex_{\vartheta_1}\exp\left\{-\frac{1}{4\sigma ^2}\sum_{j=0}^{n-1}\left[\rho
^2\,X_j^2+2\rho\,X_j\,\varepsilon _{j+1}\right]
\,\left[\1_{\left\{\BB_j\left(u_2\right)\right\}}
-\1_{\left\{\BB_j\left(u_1\right)\right\}}  \right]\right\} .
\end{align*}
Note that
$$
\1_{\left\{\BB_j\left(u_2\right)\right\}}
-\1_{\left\{\BB_j\left(u_1\right)\right\}}=
\1_{\left\{\vartheta +\frac{u_1}{n}\leq \left|X_j\right|\leq \vartheta
+\frac{u_2}{n} \right\}}\equiv\1_{\left\{\CC_j\right\}} .
$$
Further as $1-e^{-x}\leq x$,
\begin{align}
&\Ex_\vartheta
\left|Z_n^{1/2}\left(u_2\right)-Z_n^{1/2}\left(u_1\right)\right|^2\nonumber\\
&\quad =2-2\Ex_{\vartheta_1}\exp\left\{-\frac{1}{4\sigma ^2}\sum_{j=0}^{n-1}\left[\rho ^2\,X_j^2+2\rho\,X_j\,\varepsilon
_{j+1}\right]
\, \1_{\left\{\CC_j\right\}}  \right\}\nonumber \\
&\quad \leq \frac{1}{2\sigma ^2}\Ex_{\vartheta_1}\left\{\sum_{j=0}^{n-1}\left[\rho ^2\,X_j^2+2\rho\,X_j\,\varepsilon
_{j+1}\right]
\, \1_{\left\{\CC_j\right\}}  \right\}\nonumber \\
&\quad=\frac{n\rho ^2}{2\sigma
^2}\Ex_{\vartheta_1}\,X_j^2\,\1_{\left\{\CC_j\right\}}=\frac{n\rho ^2}{2\sigma
^2}\; \int_{\CC_j}^{}x^2\;f\left(\vartheta +\frac{u_1}{n},x\right)\,{\rm d}x\nonumber\\
&\quad=\frac{n\rho ^2}{2\sigma ^2} \left(\vartheta
+\frac{u_2}{n}\right)^2\,\frac{u_2-u_1}{n}\,\left[ f\left(\vartheta +\frac{u_1}{n},\tilde\vartheta_+
 \right)+f\left(\vartheta +\frac{u_1}{n},\tilde\vartheta _-\right)\right]\nonumber \\
&\quad\leq C\,\left|u_2-u_1\right|.
\label{dva}
\end{align}
We see that \eqref{2-4} is fulfilled. \hfill $\Box$

\noindent
\begin{lemma}
\label{L2} For any $p>0$ there exists a constant $C=C\left(p\right)>0$ such
 that for all values
  $u\in\left(n\left(\alpha -\vartheta\right),n\left(\beta  -\vartheta
 \right) \right)$ we have the inequality
$$
\Ex_\vartheta Z_n^{1/2}\left(u\right)\leq \frac{C}{\left|u\right|^p}.
$$
\end{lemma}
\noindent
{\bf Proof.}
We have to study the following expectation
$$
\Ex_\vartheta
Z_n^{1/2}\left(u\right)=\Ex_{\vartheta}\exp\left\{-\frac{1}{4\sigma
^2}\sum_{j=0}^{n-1}\left[\rho ^2\,X_j^2+2\rho\,X_j\,\varepsilon
_{j+1}\right] \, \1_{\left\{\BB_j\left(u\right)\right\}}  \right\}\,.
$$
Start with the  probability
\begin{align*}
&\Pb_\vartheta \left\{\ln Z_n^{1/2}\left(u\right)>-c
\left|u\right|\right\}=\Pb_\vartheta \left\{-\frac{1}{4\sigma
^2}Y_n\left(u\right)>-c\left|u\right|\right\}\\
&\quad =\Pb_\vartheta \left\{
-\frac{1}{4\sigma
^2}\sum_{j=0}^{n-1}\left[\rho ^2\,X_j^2+2\rho\,X_j\,\varepsilon
_{j+1}\right] \, \1_{\left\{\BB_j\left(u\right)\right\}}
>-c \left|u\right|\right\}\,.
\end{align*}
Write
\begin{align*}
&\Pb_\vartheta \left\{-\frac{1}{8\sigma
^2}Y_n\left(u\right)>-\frac{c}{2}\left|u\right|\right\}
 =\Pb_\vartheta \left\{ -\frac{3}{32\sigma^2}
\sum_{j=0}^{n-1}\rho ^2\,X_j^2 \, \1_{\left\{\BB_j\left(u\right)\right\}}
\right.\\
&\quad \qquad \left.
- \frac{1}{32\sigma^2}\sum_{j=0}^{n-1}\rho ^2\,X_j^2 \,
\1_{\left\{\BB_j\left(u\right)\right\}}
 - \frac{1}{4\sigma^2}\sum_{j=0}^{n-1}\rho \,X_j\,\varepsilon _{j+1} \,
\1_{\left\{\BB_j\left(u\right)\right\}}  > -\frac{c}{2}\left|u\right|
\right\} \\
&\quad \leq \Pb_\vartheta \left\{ -\frac{3}{32\sigma^2}
\sum_{j=0}^{n-1}\rho ^2\,X_j^2 \, \1_{\left\{\BB_j\left(u\right)\right\}}>-
\frac{3c}{2}\left|u\right|\right\}\\
&\qquad +\Pb_\vartheta\left\{- \frac{1}{32\sigma^2}\sum_{j=0}^{n-1}\rho ^2\,X_j^2 \,
\1_{\left\{\BB_j\left(u\right)\right\}}
 - \frac{1}{4\sigma^2}\sum_{j=0}^{n-1}\rho \,X_j\,\varepsilon _{j+1} \,
\1_{\left\{\BB_j\left(u\right)\right\}}  > c\left|u\right|
\right\} .
\end{align*}
For the last probability, by Markov inequality we have
\begin{align*}
\Pb_\vartheta\left\{- \frac{1}{32\sigma^2}\sum_{j=0}^{n-1}\rho ^2\,X_j^2 \,
\1_{\left\{\BB_j\left(u\right)\right\}}
 - \frac{1}{4\sigma^2}\sum_{j=0}^{n-1}\rho \,X_j\,\varepsilon _{j+1} \,
\1_{\left\{\BB_j\left(u\right)\right\}}  > c\left|u\right|
\right\}\leq e^{-c\left|u\right|}
\end{align*}
because
\begin{align*}
\Ex_\vartheta \exp\left\{- \frac{1}{32\sigma^2}\sum_{j=0}^{n-1}\rho ^2\,X_j^2 \,
\1_{\left\{\BB_j\left(u\right)\right\}}
 - \frac{1}{4\sigma^2}\sum_{j=0}^{n-1}\rho \,X_j\,\varepsilon _{j+1} \,
\1_{\left\{\BB_j\left(u\right)\right\}} \right\}=1.
\end{align*}
The last equality follows from the following property of the conditional
expectation
\begin{align*}
&\Ex_\vartheta \left(\left. \exp\left\{-\frac{\rho
^2}{32\sigma ^2}X_{n-1}^2\1_{\left\{\BB_{n-1}\left(u\right)\right\}}-
\frac{\rho}{4\sigma^2} \,X_{n-1}\,\varepsilon _{n} \,
\1_{\left\{\BB_{n-1}\left(u\right)\right\}}   \right\}\right| {\cal
F}_{n-1}\right)\\
&\; =\exp\left\{-\frac{\rho
^2X_{n-1}^2\1_{\left\{\BB_{n-1}\left(u\right)\right\}}}{32\sigma ^2}
\right\}\Ex_\vartheta \left(\left. \exp\left\{-
\frac{\rho X_{n-1}\1_{\left\{\BB_{n-1}\left(u\right)\right\}} }{4\sigma^2}\varepsilon _{n}
  \right\}\right| {\cal
F}_{n-1}\right)=1.
\end{align*}
Hence, it is sufficient  to study the probability
$$
\Pb_\vartheta \left\{ -
\sum_{j=0}^{n-1}X_j^2 \, \1_{\left\{\BB_j\left(u\right)\right\}}>-
c_1\,\left|u\right|\right\},
$$
where $c_1=16\,c\,\sigma ^2\,\rho ^{-2}$.

Fix some  $\kappa \in\left(0,1\right)$ and consider first the {\it local values} $u$
satisfying the condition $\left|u\right|<\,n^\kappa $. Suppose that  $u>0$
(for $u<0$ the consideration is similar). Then we have
\begin{align*}
&\Pb_\vartheta \left\{\sum_{j=0}^{n-1} \left[\vartheta ^2-X_j^2\right]
\1_{\left\{\BB_j\left(u\right)
\right\}}   -\vartheta ^2\sum_{j=0}^{n-1}  \1_{\left\{\BB_j\left(u\right)
\right\}}      \geq -c_1\,u\right\}\\
&\quad =\Pb_\vartheta \left\{\sum_{j=0}^{n-1} \left[X_j^2-\vartheta ^2\right]
\1_{\left\{\BB_j\left(u\right)
\right\}}   +\vartheta ^2\sum_{j=0}^{n-1}  \1_{\left\{\BB_j\left(u\right)
\right\}}      \leq c_1\,u\right\}\\
&\quad \leq \Pb_\vartheta \left\{\sum_{j=0}^{n-1}  \1_{\left\{\BB_j\left(u\right)
\right\}}      \leq c_2\,u\right\},\qquad c_2=\frac{c_1}{\vartheta ^2},
\end{align*}
because $\left[X_j^2-\vartheta ^2\right]
\1_{\left\{\BB_j\left(u\right)
\right\}}\geq 0$. Recall that the sum $\sum_{j=0}^{n-1}  \1_{\left\{\BB_j\left(u\right)
\right\}}    $ converges to the Poisson process of intensity $\lambda
=2f\left(\vartheta ,\vartheta \right)$, hence the last probability has to be
small for the values $c_2<\lambda $. Let $Y_j=\1_{\left\{\BB_j\left(u\right)
\right\}}  -\Ex_\vartheta \1_{\left\{\BB_j\left(u\right)
\right\}}  $ and note that
\begin{align*}
\Ex_\vartheta \1_{\left\{\BB_j\left(u\right)\right\}} & =\int_{\vartheta \leq
\left|x\right|\leq \vartheta +\frac{u}{n}}^{} f\left(\vartheta ,x\right)\,{\rm
d}x =\frac{u}{n}\left[f\left(\vartheta ,-\tilde\vartheta_1
\right)+f\left(\vartheta ,\tilde\vartheta_2 \right)\right] \\
&=2\,\frac{u}{n}\,f\left(\vartheta ,\vartheta
\right)\left(1+o\left(1\right)\right) \geq \frac{u}{n}\,f\left(\vartheta ,\vartheta
\right)\,,
\end{align*}
where the last inequality holds for $n\geq n_1$. Recall that the function
$f\left(\vartheta ,x\right)$ is even and $f\left(\vartheta ,-\vartheta
\right)=f\left(\vartheta ,\vartheta\right) $.  Further,
\begin{align*}
\Pb_\vartheta \left\{\sum_{j=0}^{n-1}  \1_{\left\{\BB_j\left(u\right)
\right\}}      \leq c_2\,u\right\}&\leq \Pb_\vartheta
\left\{\left|\sum_{j=0}^{n-1}  Y_j  \right| \geq \left(f\left(\vartheta
,\vartheta\right)- c_2\right) \,u\right\}\\
&=\Pb_\vartheta
\left\{\left|\sum_{j=0}^{n-1}  Y_j  \right| \geq c_3 \,u\right\}\leq
\frac{\Ex_\vartheta  \left|\sum_{j=0}^{n-1}  Y_j  \right|^{2p}}{c_3^{2p} \,u^{2p} },
\end{align*}
where we chose such $c$ that $c_3=f\left(\vartheta ,\vartheta \right)-16c\sigma ^2\rho
^{-2}\vartheta ^{-2}>0 $.  To estimate the last expectation we apply the
inequality of Dedeker and Doukhan (see (8.1) in \cite{DD03}):
\begin{equation}
\label{DD}
\Ex_\vartheta  \left|\sum_{j=0}^{n-1}  Y_j  \right|^{2p}\leq
\left(4\,p\,n\sum_{j=0}^{n-1} \left[\Ex_\vartheta \left|Y_0\,\left.\Ex_\vartheta
\left(Y_j\right|{\cal F}_0\right)\right|^p\right]^{\frac{1}{p}}\right)^p.
\end{equation}
Write
\begin{align*}
\left.\Ex_\vartheta
\left(Y_j\right|{\cal F}_0\right)=\left.\Ex_\vartheta
\left( \left.\Ex_\vartheta  \left(Y_j\right|{\cal F}_{j-1} \right)\right|{\cal F}_0\right)
\end{align*}
and let $g\left(x\right)=\rho
_1\,x\1_{\left\{\left|x\right|<\vartheta \right\}}+ \rho
_2\,x\1_{\left\{\left|x\right|\geq \vartheta \right\}}$, then
\begin{align*}
\left.\Ex_\vartheta  \left(Y_j\right|{\cal F}_{j-1}
\right)&=\left.\Ex_\vartheta  \left(
\1_{\left\{\BB_{j}\left(u\right)\right\}}        \right|{\cal F}_{j-1}
\right) -\Ex_\vartheta   \1_{\left\{\BB_{j}\left(u\right)\right\}} \\
&=    \left.\Ex_\vartheta  \left(
\1_{\left\{\vartheta \leq \left|g\left(X_{j-1}\right)+\varepsilon _j
\right|\leq \vartheta +\frac{u}{n}\right\}}        \right|{\cal F}_{j-1}
\right) -\Ex_\vartheta   \1_{\left\{\BB_{j}\left(u\right)\right\}}\\
&=\int_{\vartheta \leq \left|g\left(X_{j-1}\right)+x  \right|\leq \vartheta
+\frac{u}{n}}^{} \varphi \left(x\right)\,{\rm d}x +\int_{\vartheta \leq \left|x
\right|\leq \vartheta
+\frac{u}{n}}^{} f \left(\vartheta ,x\right)\,{\rm d}x\\
&=\frac{u}{n}\; \left[\varphi \left(\vartheta
-g\left(X_{j-1}\right)+\frac{\tilde u}{n}\right)+ \varphi \left(\vartheta
+g\left(X_{j-1}\right)-\frac{\tilde{\tilde u}}{n}\right)\right]\\
&\quad +\frac{u}{n}\; \left[f \left(\vartheta ,\vartheta
+\frac{\bar u}{n}\right)+ f \left(\vartheta ,\vartheta
-\frac{\bar{\bar u}}{n}\right)\right]=\frac{u}{n}\,A\left(X_{j-1}\right),
\end{align*}
where $\varphi \left(\cdot \right)$ is the density function of the Gaussian
r.v. $\varepsilon _j$, i.e., $\varphi \left(\cdot \right)\sim {\cal
N}\left(0,\sigma ^2\right)$) and $\bar u,\tilde u, \tilde{\tilde u}$ are some
values between $0$ and $u$. Note that the function $A\left(x\right)$ defined
by the last equality is bounded and $\Ex_\vartheta A\left(X_{j-1}\right)=0$.

Hence
\begin{align*}
\left.\Ex_\vartheta \left(Y_j\right|{\cal
F}_0\right)=\frac{u}{n}\,\left.\Ex_\vartheta \left(A\left(X_{j-1}\right)
\right|{\cal F}_0\right)
\end{align*}
and
\begin{align*}
&\sum_{j=0}^{n-1} \left[\Ex_\vartheta \left|Y_0\,\left.\Ex_\vartheta
\left(Y_j\right|{\cal F}_0\right)\right|^p\right]^{\frac{1}{p}}=\frac{u}{n}\,
\sum_{j=0}^{n-1} \left[\Ex_\vartheta \left|Y_0\,\left.\Ex_\vartheta
\left(A\left(X_{j-1}\right)\right|{\cal
F}_0\right)\right|^p\right]^{\frac{1}{p}}\\
&\qquad \leq C\;\frac{u}{n}\, \sum_{j=0}^{n-1}\alpha
\left(j-1\right)=C\;\frac{u}{n}\, \sum_{j=0}^{n-1}\gamma ^{j-1 }\leq C\;\frac{u}{n}\,,
\end{align*}
where we used the geometrical ergodicity of $\left(X_j\right)_{j\geq 1}$:
$\alpha \left(j\right) \leq \gamma ^j$,  $0< \gamma <1$ (see \cite{CT91}) and  the
inequality of Ibragimov
\begin{align*}
&\left\|\Ex_\vartheta \left(\left. A\left(X_{j}\right)  \right| {\cal
F}_0\right)-\Ex_\vartheta A\left(X_{j} \right)  \right\|_p\\
&\qquad \qquad \leq
C\left\|\left.\Ex_\vartheta \left(A\left(X_{j}\right)\right| {\cal
F}_0\right)-\Ex_\vartheta A\left(X_{j} \right)  \right\|_1^{1/p}
\leq C\,\alpha
\left(j\right)^{1/p}
\end{align*}
(see Bradley \cite{Br1}, Theorem 4.4, (a2)).

Finally, we obtain  (for $\left|u\right|\leq n^\kappa $) the estimate
\begin{align*}
\Pb_\vartheta \left\{\sum_{j=0}^{n-1}  \1_{\left\{\BB_j\left(u\right)
\right\}} \leq c_2\left|u\right|  \right\}    \leq \frac{C}{\left|u\right|^p}\,.
\end{align*}

Consider now the case $\left|u\right|>n^{\kappa }$. Of course,
$\left|u\right|\leq \left(\beta -\alpha \right)n$.
We have
\begin{align}
&\Pb_\vartheta
\left\{-\sum_{j=0}^{n-1}X_j^2\1_{\left\{\BB_j\left(u\right)
\right\}}>-c_1\left|u\right|\right\}\nonumber\\
&\qquad \leq \Pb_\vartheta \left\{\left|\sum_{j=0}^{n-1}\left(X_j^2\,\1_{\left\{\BB_j\left(u\right)
\right\}}-\Ex_\vartheta X_j^2\,\1_{\left\{\BB_j\left(u\right) \right\}}
\right)\right|\geq c_1\left|u\right| \right\},
\label{2-8}
\end{align}
because
\begin{align*}
\Ex_\vartheta X_j^2\,\1_{\left\{\BB_j\left(u\right)
\right\}}=\int_{\vartheta
<\left|x\right|<\vartheta +\frac{\left|u\right|}{n}}^{}x^2f\left(\vartheta ,x\right)\,{\rm
d}x \geq \alpha  ^2\;\frac{\left|u\right|}{n} \inf_{\alpha \leq x\leq \beta
}f\left(\vartheta ,x\right)\geq \frac{c_4}{n}\,\left|u\right|,
\end{align*}
and the constant $c$ is chosen such that $c_4>2c_1$.
The last  probability in \eqref{2-8} can be estimated
with the help of the following lemma.

\begin{lemma} {\rm (Rosenthal's moment inequality)}
\label{Ros}
 Let $\left(Z_j\right)_{j\geq 1}$ be zero mean mixing series satisfying
the condition:
there exist $\varepsilon >0$
and $c\in 2\NN$, $c>2p>2$, such that
\begin{equation}
\label{mix1}
\sum_{r=1}^{\infty }\left(r+1\right)^{c-2}\left[\alpha
\left(r\right)\right]^\frac{\varepsilon }{c+\varepsilon }<\infty  ,
\end{equation}
where $\alpha \left(r\right)$ is the $\alpha $-mixing coefficient, then
\begin{equation}
\label{mix2}
\Ex\left|\sum_{j=1}^{n}Z_j\right|^{2p}\leq
C\left[n\;\left(\Ex\left|Z_1\right|^{2p+\varepsilon
}\right)^{\frac{2p}{2p+\varepsilon }} + n^{p}\left(\Ex Z_1^{2+\varepsilon
}\right)^\frac{2p }{2+\varepsilon } \right] .
\end{equation}
 \end{lemma}
For the proof  see \cite{PD}, p.26.

As the process $X_j$ is geometrically mixing (see \cite{FY03}, Theorem 2.4), hence
condition \eqref{mix1} is fulfilled with any $c>0$ and $\varepsilon >0$. We apply
 \eqref{mix2}  with
$$
Z_j=X_j^2\,\1_{\left\{\BB_j\left(u\right)
\right\}}-\Ex_\vartheta X_j^2\,\1_{\left\{\BB_j\left(u\right) \right\}}
.
$$
 Obviously, $e_n\left(u\right)\equiv\Ex_\vartheta
 X_j^2\,\1_{\left\{\BB_j\left(u\right)\right\}}  \leq \beta ^2$ and
 $e_n\left(u\right)\leq C\,\frac{\left|u\right|}{n}$. We suppose for simplicity
 that $u>0$,
\begin{align*}
\Ex_\vartheta \left|Z_1\right|^{2p+\varepsilon}&=\int_{\vartheta
<\left|x\right|<\vartheta +\frac{u}{n}}^{} \left|x^2-
e_n\left(u\right) \right|^{2p+\varepsilon}f\left(\vartheta ,x\right) {\rm
d}x\\
&\quad
+e_n\left(u\right)^{2p+\varepsilon}\left(1-\int_{\vartheta
<\left|x\right|<\vartheta +\frac{u}{n}}^{}f\left(\vartheta ,x\right) {\rm
d}x  \right)\\
&\leq
C_1\;\frac{u}{n}+C_2\left|\frac{u}{n}\right|^{2p+\varepsilon }
\end{align*}
and similarly
$$
\Ex_\vartheta \left|Z_1\right|^{2+\varepsilon}\leq
C_3\;\frac{u}{n}+C_4\left|\frac{u}{n}\right|^{2+\varepsilon }.
$$
Hence,
\begin{align*}
&\Ex_\vartheta\left|\sum_{j=0}^{n-1}\left(X_j^2\,\1_{\left\{\BB_j\left(u\right) \right\}}-\Ex_\vartheta
X_j^2\,\1_{\left\{\BB_j\left(u\right) \right\}} \right)\right|^{2p}\\
&\qquad\qquad \leq
C\,\left(n\,\left|\frac{u}{n}\right|^{\frac{2p}{2p+\varepsilon }} +
n\left|\frac{u}{n}\right|^{2p}+n^p\left|
\frac{u}{n}\right|^{\frac{2p}{2+\varepsilon }}
+n^p\left|\frac{u}{n}\right|^{2p} \right)\\
 &\qquad\qquad
C\,\left(n^{\frac{\varepsilon }{2p+\varepsilon
}}\left|u\right|^{1-\frac{\varepsilon }{2p+\varepsilon
}}+n^{1-2p}\left|u\right|^{2p} +n^{\frac{p\varepsilon }{2+\varepsilon
}}\left|u\right|^{p-\frac{\varepsilon
}{2+\varepsilon}}+\frac{\left|u\right|^{2p}}{n^p} \right)\\
&\qquad\qquad\leq
C\,\left( \left|u\right|^{\frac{\varepsilon }{\kappa \left(2p+\varepsilon\right)
}+1-\frac{\varepsilon }{2p+\varepsilon
}}+\left|u\right|^{\frac{1-2p}{\kappa }+2p}+\left|u\right|^{\frac{p\varepsilon }{\kappa \left(2+\varepsilon\right)
}+p-\frac{p\varepsilon }{2+\varepsilon}}+\left|u\right|^{p} \right)\leq
C\,\left|u\right|^{p},
\end{align*}
where we have used the relations $n<\left|u\right|^{1/\kappa }$,
$\left|u\right|\leq \left(\beta -\alpha \right)n$ and have chosen sufficiently
small $\varepsilon $ (or sufficiently large $p$).

By Chebyshev inequality
\begin{align*}
&\Pb_\vartheta \left\{\left|\sum_{j=0}^{n-1}\left(X_j^2\,\1_{\left\{\BB_j\left(u\right)
\right\}}-\Ex_\vartheta X_j^2\,\1_{\left\{\BB_j\left(u\right) \right\}}
\right)\right|\geq {c_1}\left|u\right| \right\}\leq
\frac{C}{\left|u\right|^{p}}.
\end{align*}

The estimates obtained above allow us to write the following expression: for
any $p>1$ and all $u\in \left(n\left(\alpha -\vartheta\right) ,n\left(\beta
-\vartheta\right)\right)$, there exist constants $c>0$ and $C>0$ such that
\begin{equation}
\label{ld}
\Pb_\vartheta \left\{Z_n\left(u\right)>e^{-c\left|u\right|}\right\}\leq \frac{C}{\left|u\right|^p}.
\end{equation}
For the expectation, note that
\begin{align*}
\Ex_\vartheta Z_n^{1/2}\left(u\right)&=\Ex_\vartheta
Z_n^{1/2}\left(u\right)\1_{\left\{Z_n^{1/2}\left(u\right)\geq
e^{-\frac{c}{2}\left|u\right|}\right\}}+\Ex_\vartheta
Z_n^{1/2}\left(u\right)\1_{\left\{Z_n^{1/2}\left(u\right)<
e^{-\frac{c}{2}\left|u\right|}\right\}}\\
&\leq \left(\Ex_\vartheta
Z_n\left(u\right)  \Pb_\vartheta
\left\{Z_n\left(u\right)>e^{-c\left|u\right|}\right\}\right)^{1/2}
+e^{-\frac{c}{2}\left|u\right|}\leq  \frac{C}{\left|u\right|^{p/2}}.
\end{align*}
Recall that this estimate is valid for any $p>1$, hence \eqref{2-5} is
verified. Therefore the required conditions are fulfilled and the Bayesian
estimate satisfied all of the properties stipulated in Theorem 1 (see Theorem 1.10.2,
\cite{IH81}). \hfill $\Box$

\section[short title]{Simulations}

We obtain the density functions of limit distributions of the MLE and
Bayesian estimators by the following simulations. The limit likelihood ratio
is
$$
Z\left(u\right)=\exp\left\{-\frac{\rho ^2\vartheta ^2}{2\sigma ^2}\;
N_+\left(u\right) -\frac{\rho \vartheta }{\sigma
^2}\,\sum_{l=0}^{N_+\left(u\right)} \varepsilon_l^+  \right\}.
$$
for $u\geq 0$ and
$$
Z\left(u\right)=\exp\left\{-\frac{\rho ^2\vartheta ^2}{2\sigma ^2}\;
N_-\left(-u\right) -\frac{\rho \vartheta }{\sigma
^2}\,\sum_{l=0}^{N_-\left(-u\right)} \varepsilon_l^-  \right\}.
$$
for $u\leq  0$. Here $N_+\left(\cdot \right)$ and $N_-\left(\cdot \right)$ are
independent Poisson processes of intensity
$
\lambda=2f\left(\vartheta ,\vartheta \right)
$
and the Gaussian random variables $\varepsilon _l^+ $, $\varepsilon _l^- \sim {\cal
N}\left(0,\sigma ^2\right)$ are independent, $\varepsilon _0^\pm=0$.

Denote
$$
\gamma =\frac{\left(\rho _2-\rho _1\right)\vartheta }{\sigma }, \qquad
\epsilon_l^\pm=\frac{-\varepsilon _l^\pm}{\sigma }  ,\qquad u=\frac{v}{\lambda
},\quad \nu _\pm\left(v\right) =N_\pm\left(\frac{v}{\lambda }\right).
$$
Then the Poisson processes $\nu _+\left(v\right),v\geq 0$ and $\nu
_-\left(v\right),\geq 0$ have intensity 1 and the limit likelihood ratio
\begin{align*}
Z_\gamma \left(v\right)=\begin{cases}
\exp\left\{\gamma \,\sum_{l=0}^{\nu _{-}\left(-v\right)} \epsilon_l^{-} -\frac{\gamma ^2}{2}
\nu _-\left(-v\right) \right\}, &\text{if $v\leq 0$,}\\
\exp\left\{ \gamma \,\sum_{l=0}^{\nu _+\left(v\right)}\,\, \epsilon_l^+ -\,\,\frac{\gamma ^2}{2}
\nu _+\left(v\right)
\right\}, &\text{if $v> 0.$}
			\end{cases}
\end{align*}
Now the limit process $Z_\gamma \left(v\right)$ only depends on one parameter
($\gamma $) and the limit random variables $\hat u$ and $\tilde u$ can be
written as
$$
\hat u=\frac{\hat u_\gamma }{\lambda },\quad \tilde u=\frac{\tilde u_\gamma
}{\lambda },\qquad  \tilde u_\gamma=\frac{\int_{-\infty }^{\infty }v\,Z_\gamma
\left(v\right)\,{\rm d}v  }{\int_{-\infty }^{\infty }Z_\gamma
\left(v\right)\,{\rm d}v}\,,
$$
in obvious notation.

The next problem is to find the function $f\left(\vartheta ,x\right)$, where
$f\left(\vartheta ,x\right)$ is the stationary density function of $X_j$. As
$$
X_{j+1}=\rho_1
\,X_j+\rho
\,X_j\,\1_{\left\{\left|X_j\right|>\vartheta \right\}}+\varepsilon _{j+1}
 $$
where $X_j$ and $\varepsilon _{j+1}$ are independent, we obtain the
convolution equation
$$
f\left(\vartheta ,y\right)=\frac{1}{\sqrt{2\pi \sigma ^2}}\int_{-\infty }^{\infty
}g\left(\vartheta ,x\right) e^{-1/2\sigma ^2\left(y-x\right)^2}{\rm d}x\,.
$$
Herein, we denote the density function of
$\rho_1 \,X_j+\rho \,X_j\,\1_{\left\{\left|X_j\right|>\vartheta \right\}}$
by $g\left(\vartheta ,x\right)$. This density can be
expressed as a function of $f\left(\cdot \right)$, which is a solution to a
corresponding integral equation, see also Chan and Tong \cite{CT86}.  Specifically,
observe that
$$
g\left(\vartheta ,x\right)=\frac{1}{\rho_1 }f\left(\vartheta ,\frac{x}{\rho_1
}\right)\;\1_{\left\{\frac{\left|x\right|}{\rho _1}<\vartheta \right\}}+\frac{1}{\rho_2
}f\left(\vartheta ,\frac{x}{\rho_2
}\right)\;\1_{\left\{\frac{\left|x\right|}{\rho _2}\geq \vartheta \right\}}.
$$
Hence, the integral equation is
\begin{align}
&f\left(\vartheta ,y\right)=\int_{-\infty }^{\infty }\left[\frac{1}{\rho_1
}f\left(\vartheta ,\frac{x}{\rho_1
}\right)\;\1_{\left\{\frac{\left|x\right|}{\rho _1}<\vartheta
\right\}}\right.\nonumber\\
&\qquad \qquad \quad \left.+\frac{1}{\rho_2
}f\left(\vartheta ,\frac{x}{\rho_2
}\right)\;\1_{\left\{\frac{\left|x\right|}{\rho _2}\geq \vartheta
\right\}}\right] \varphi \left(y-x\right){\rm d}x\nonumber\\
&\qquad =\int_{-\infty }^{\infty } f\left(\vartheta ,x\right) \left[\varphi
\left(y-x\rho _1\right)\;\1_{\left\{\left|x\right|<\vartheta
\right\}} + \varphi
\left(y-x\rho _2\right)\;\1_{\left\{\left|x\right|\geq \vartheta
\right\}} \right] {\rm d}x .
\label{dens}
\end{align}

Solution to this equation at the point $\vartheta $ is the intensity
$\lambda =2f\left(\vartheta ,\vartheta \right)$ of the corresponding Poisson
processes.  Therefore the value $f\left(\vartheta ,\vartheta \right)$ satisfies
the integral equation
$$
f\left(\vartheta ,\vartheta \right)=\int_{-\infty }^{\infty } f\left(\vartheta ,x\right) \left[\varphi
\left(\vartheta -x\rho _1\right)\;\1_{\left\{\left|x\right|<\vartheta
\right\}} + \varphi
\left(\vartheta -x\rho _2\right)\;\1_{\left\{\left|x\right|\geq \vartheta
\right\}} \right] {\rm d}x\,,
$$
where  $\varphi \left(\cdot
\right)$ is Gaussian ${\cal N}\left(0,\sigma ^2\right)$ density. We see that
$f\left(\vartheta ,\vartheta \right)>0 $.  To visualize the properties of
the sample path $Z(u)$ and the invariant density $f(\theta,\vartheta )$, we conduct a
simulation experiment by taking
$\varepsilon^{\pm}_l$ to be i.i.d. standard normal random variables.  The parameters
used are $\vartheta=2, \rho_1=0.15, \rho_2=0.95, \sigma=1$ and $\lambda=0.5$.
A Gaussian kernel is used to estimate the form of $f(\theta,2 )$ based on $50,000$
observations of $X_t$ generated from model~\eqref{2-6} with $\varepsilon_t$
being i.i.d. ${\cal N}\left(0,1\right)$ random variables.  The plots of $Z(u)$
and $f(\theta,2)$ are given in Figures 1 and 2 respectively.

\begin{center}
{\em Figure 1.}\\
{\em Figure 2.}
\end{center}

For the maximum likelihood estimate, note that the maximum values
of $Z\left(u\right)$ form an interval [$\hat u_m,\hat u_M$] with
length $ \left|\hat u_M-\hat u_m\right|=\eta $, where
$\eta $ is an exponential random variable with probability density
$2f_\vartheta e^{-2f_\vartheta\,x},x\geq 0$,  $f_\vartheta
=f\left(\vartheta ,\vartheta \right)$. We can take any value of $u$ from this
interval, the middle point $\hat u=\frac{\hat u_M+\hat u_m\ }{2}$, say.
To have its density function we only need to simulate the exponential and the Gaussian
independent random variables which will generate $\hat u_1,\ldots
\hat u_N$.  The historgram of $\hat{u}$ based on $20,000$ simulated values
of $\hat{u}$ is plotted in Figure 3.  As can be seen clearly, the MLE performs
reasonably well and converges to zero very fast.  The sample mean of the simulated
$\hat{u}$ is $0.01$ with a standard deviation $0.033$.

For the Bayesian estimators we first calculate the integral
\begin{align*}
J_+=\int_{0}^{\infty }u\,Z\left(u\right)\,{\rm d}u&=\sum_{l=0}^{\infty
}\int_{u_l}^{u_{l+1}}u\;e^{-\frac{\rho ^2\vartheta ^2}{2\sigma
^2}\;l-\frac{\rho \vartheta }{\sigma
^2}\,\sum_{r=0}^{l} \varepsilon_r}{\rm d}u\\
&=\frac{1}{2}\sum_{l=0}^{\infty
}e^{-\frac{\rho ^2\vartheta ^2}{2\sigma
^2}\;l-\frac{\rho \vartheta }{\sigma
^2}\,\sum_{r=0}^{l} \varepsilon_r} \left(u_{l+1}^2- u_l^2\right).
\end{align*}
Here $\varepsilon _r\sim {\cal N}\left(0,\sigma ^2\right)$ and
$u_l=\sum_{r=0}^{l}\eta _r$. By a similar way we have
$$
I_+=\int_{0}^{\infty }Z\left(u\right)\,{\rm d}u=\sum_{l=0}^{\infty
}e^{-\frac{\rho ^2\vartheta ^2}{2\sigma
^2}\;l-\frac{\rho \vartheta }{\sigma
^2}\,\sum_{r=0}^{l} \varepsilon_r} \left(u_{l+1}- u_l\right).
$$
The limit random variable is
$$
\tilde u=\frac{J_-+J_+}{I_-+I_+}\,,
$$
with obvious notation.
To understand the behavior of the Bayesian estimator, we simulate the Bayesian estimator for
$20,000$ times with the histogram of $\tilde{u}$ given in Figure 4.  From this figure, it
is clearly seen that the Bayesian estimator converges to the expected value zero.  The sample
mean is $-0.0026$ with a standard deviation $0.028$.   It is interesting to see that this
simulation results are consistent with the theory that the limit variances of the MLE and BE
$$
d_{MLE}^2 \sim \frac{1}{N}\sum_{q=1}^{N}\hat u_q^2 \ \mbox{ and } \ d_{BE}^2 \sim
\frac{1}{N}\sum_{q=1}^{N}\tilde u_q^2
$$
satisfy
$$
d_{MLE}^2 >d_{BE}^2 .
$$
Note that it follows from the symmetry of the limit process, the random
variables $\hat u$ and $\tilde u$ satisfy $\Ex_\vartheta \hat
u=0=\Ex_\vartheta \tilde u$.

For 20,000 simulated estimators, we obtain the limit variances as
$d_{MLE}^2= 22.83\pm 0.68$ and $d_{BE}^2=16.79\pm 0.39$.
These values concur  with the theoretical results that the
Bayesian estimator outperforms the MLE.

\begin{center}
{\em Figure 3.}\\
{\em Figure 4.}
\end{center}

To examine the finite sample performance of the test statistics, we computed
the critical values of the limit distributions based on the MLE and the BE using
the same set of parameters as given in Figures 3.
The sizes are chosen for commonly used test statistics and the limiting
values are given in the first two rows of Table 1.  As can be seen, both the MLE
and BE procedures perform reasonably well and are in close agreement.
Furthermore, the numbers in the last row of Table 1
are the critical values computed from the test statistics in
(\ref{2-2}), which are directly simulated from model (\ref{2-1}) using the same set
of parameters.  It is seen that the critical values generated from the simulated
statistics agree remarkably well with the critical values computed from the limit
distributions in Table 1 based on MLE.  In summary, Table 1 demonstrates
the usefulness of the limit distributions in computing the critical values.  If
one needs to conduct a test for another set of parameters, then a similar table
can be computed and the programming code is available from the authors upon request.

\begin{center}
\begin{table}
\begin{tabular}[c]{|ccccccccc|}
\hline
 & 0.025 & 0.05 & 0.075 & 0.1 & 0.90 & 0.925 & 0.95 & 0.975\\
\hline MLE & $-9.66$ & $-6.64$ & $-5.28$ & $-4.46$ &
 4.46 & 5.38& 6.87 & 9.84\\
BE & $-8.44$ & $-6.29$ & $-5.07$ & $-4.27$ & 4.21 &
5.09 & 6.26 & 8.43\\
Simulated Values & $-9.70$ & $-6.88$ & $-5.48$ & $-4.64$ &
4.90 & 5.85 & 7.55 & 10.28 \\
\hline
\end{tabular}
\caption{Critical values for $\lambda=0.5$, $\theta=2$, $\sigma=1$,
$\rho_1=0.95$ and $\rho_2=0.15$.}
\end{table}
\end{center}

\section[short title]{One-sided threshold} In the nonlinear time series
literature, the threshold AR model usually takes the form (see, e.g.,
\cite{Ch93}, \cite{Han97} and \cite{To90})
\begin{equation}
\label{4-1}
X_{j+1}=\rho _1X_j\1_{\left\{X_j<\vartheta \right\}}+\rho
_2X_j\1_{\left\{X_j\geq \vartheta \right\}} \varepsilon _{j+1}
\end{equation}
Note that the study of these one-sided threshold models is no more
complicated than \eqref{2-1} because the log-likelihood ratio
$$
\ln Z_n\left(u\right)=\frac{L\left(\vartheta
+\frac{u}{n},X^n\right)}{L\left(\vartheta ,X^n\right)} =-\frac{1}{2\sigma
^2}\;Y_n\left(u\right)
$$
 depends on the stochastic process (for fixed $u>0$)
$$
Y_n\left(u\right)=\sum_{j=0}^{n-1}\left(\rho ^2X_j^2+2\rho X_j
\varepsilon _{j+1}\right)\;\1_{\left\{\vartheta <X_j<\vartheta +\frac{u}{n}\right\}}
$$
(see \eqref{lr}), which is approximated by the process
$$
Y_n^\circ\left(u\right)=\sum_{j=0}^{n-1}\left(\rho ^2\vartheta ^2+2\rho \vartheta
\varepsilon _{j+1}\right)\;\1_{\left\{\vartheta <X_j<\vartheta +\frac{u}{n}\right\}}.
$$
Here we use the same notations as before and add the condition that $\left|\rho
_1\right|<1$.  Comparison with \eqref{lro} shows that the
factor $\sgn\left(X_j\right)$ no longer exists and this simplifies matters
much in the application of the limit theorems.  Specifically, similar to \eqref{ll},
the corresponding limit for $Y_n^\circ\left(u\right)$ becomes
$$
Y_n^\circ\left(u\right)\Longrightarrow
Y_+\left(u\right)=\sum_{l=0}^{N_+\left(u\right)}\left[ \rho ^2\vartheta ^2+2\rho \vartheta \varepsilon _l^+\right],
$$
the only difference is: instead of $2f\left(\vartheta ,\theta \right)$, the intensity of the Poisson process
$N_+\left(u\right)$ is  $\lambda_+=f\left(\vartheta ,\vartheta \right)$.

The inequalities for the process $Z_n\left(u\right)$ obtained in Lemmas~\ref{L1} and \ref{L2}
can be obtained for $Z_n\left(u\right)$ of the process
\eqref{4-1} exactly the same way as in this paper.
Consequently, the asymptotic behavior of the Bayesian estimator $\tilde\vartheta _n$
for model \eqref{4-1} is the same as that described in the
Theorem \ref{T1} with the slightly difference due to the form of the limit
likelihood ratio $Z\left(u\right)$, where the intensity of the Poisson process is
now $f\left(\vartheta ,\vartheta \right)$, not $2f\left(\vartheta ,\theta \right)$.

\section[short title]{Discussion}
Let us explain heuristically why the choice \eqref{choice} for the MLE
is better than other types of the form $\hat u_\gamma =\gamma u_m+\left(1-\gamma \right)u_M$
with $\gamma \not = 1/2, \gamma \in\left[0,1\right]$. The interval
$\left[u_m,u_M\right]$ can be on the positive, negative parts of $R$ or it can
be $\left[u_1^-,u_1^+\right]$, where $u_1^-$ and $u_1^+$ are the first event
of the Poisson processes $N_-\left(\cdot \right)$ and $N_+\left(\cdot \right)$
respectively. If $\left[u_m,u_M\right]=\left[u_1^-,u_1^+\right]$, then the random
variables $-u_1^-=\zeta _1$ and $u_1^+=\zeta _2$ are independent exponential
with the  parameter $\lambda =2f\left(\vartheta ,\vartheta \right)$.

If this interval is on the negative part, then $\left[u_{\hat l+1}^-
,u_{\hat l}^-\right]$ ($u_l^-$ is the $l$-th event of the Poisson process
$N_-\left(\cdot \right)$) has  random length $\hat l$ and
$$
\hat u_\gamma =\gamma u_{\hat
l+1}^-+\left(1-\gamma \right)u_{\hat l}^-=\gamma \left(u_{\hat
l}^--\zeta\right) +\left(1-\gamma \right)u_{\hat l}^-
=u_{\hat l}^--\gamma \zeta ,
$$
where $\zeta $ is exponential random variable with parameter $\lambda $. For the positive part
$$
\hat u_\gamma =\gamma u_{\hat
l}^++\left(1-\gamma \right)u_{\hat l+1}^+=\gamma u_{\hat
l}^+ +\left(1-\gamma \right)\left(u_{\hat l}^++\zeta \right) =u_{\hat
l}^++\left(1-\gamma\right) \zeta .
$$
Denote $p_0$ to be the probability that the maximum of the random process
$Z\left(\cdot \right)$ is on the interval $\left[u_1^-,u_1^+ \right]$. The
positive and negative intervals are equiprobable, hence their probabilities
$p_-, p_+$ satisfy the relations $p_-= p_+=\left(1-p_0\right)/2\equiv p$. We then write
\begin{align*}
\Ex_\vartheta \hat u_\gamma^2&=p_- \,\Ex_\vartheta \left(u_{\hat l}^--\gamma
\zeta\right)^2+
p_0 \,\Ex_\vartheta\left(\gamma \zeta _1-\left(1-\gamma \right)\zeta _2\right)^2\\
&\quad
+p_+ \,\Ex_\vartheta \left(u_{\hat l}^++\left(1  -\gamma\right) \zeta\right)^2=
2p \,\Ex_\vartheta \left(u_{\hat l}^+\right)^2+ 2p \, \Ex_\vartheta \left(u_{\hat
l}^+\right)\;\Ex_\vartheta \zeta \\
&\quad
+\frac{2 \,p}{\lambda ^2} \,\left(2\gamma ^2-2\gamma +1\right)+\frac{2\,p_0}{\lambda ^2} \,\left(4\gamma
^2-4\gamma +1\right)
\end{align*}
and direct calculations show
$$
\min_{\gamma \in \left[0,1\right]}\Ex_\vartheta \left(\hat u_\gamma \right)^2
=\Ex_\vartheta \left(\hat u \right)^2 .
$$

Note that in this problem of parameter estimation, it is possible to introduce
the notion of asymptotic efficiency of estimators. The lower bound on the risk
of all estimators $\bar\vartheta _n$ for the quadratic loss function is as follows:
\begin{equation}
\label{lbd}
\Liminf_{\delta \rightarrow 0}\Liminf_{n \rightarrow\infty
}\sup_{\left|\vartheta -\vartheta _0\right|<\delta }n^2 \,\Ex_\vartheta
\left(\bar \vartheta _n-\vartheta \right)^2\geq \Ex_{\vartheta _0}\left(\tilde
u\right)^2.
\end{equation}

This bound  follows from the results of the Section 1.9
in \cite{IH81}. We just note that the second moment
$\Ex_{\vartheta}\left(\tilde u\right)^2$ is a continuous function of $\vartheta
$.

As $\Ex_{\vartheta}\left(\tilde u\right)^2$ is the limit of the Bayesian
estimator, we can think of these estimators having smaller limit error than the
MLE (as in singular estimation problems). To prove this asymptotic
efficiency of the Bayesian estimators, we need to show that the convergence of
the second moments is uniform in $\vartheta $ on compacts. The corresponding
uniform  estimates
on the process $Z_n\left(\cdot \right)$ can be easily verified and what remains
to be done is to
establish the uniform version of the convergence of finite dimensional
distributions, which can also be verified.

Another possible generalization is to consider the case where $\{\varepsilon _j\}_{j\geq 1}$ are
independent random variables
with a known density function satisfying  some regularity conditions. Then the
estimator $\hat \vartheta _n$ (defined by \eqref{mle}) becomes the least squares
estimator and $\tilde\vartheta_n$ (defined by \eqref{be}) is no longer Bayesian,
but becomes another estimator having desirable asymptotic properties.  The behaviour of these
estimators can be similarly studied and their limit distributions can
be defined via the corresponding limit process $Z\left(\cdot \right)$ when
$\varepsilon _l^\pm$ are no longer Gaussian.

Note that a continuous-time analogue of TAR model is prescribed by the following
stochastic differential equation
\begin{equation*}
{\rm d}X_t=-\varrho_1 \,X_t\,\1_{\left\{\left|X_t\right|<\vartheta \right\}}{\rm
d}t-\varrho_2 \,X_t\,\1_{\left\{\left|X_t\right|\geq \vartheta \right\}}{\rm
d}t+\sigma  \,{\rm d}W_t,\quad  0\leq t\leq T,
\end{equation*}
where $\varrho_1\not =\varrho_2 >0$.
This model can be called {\it Threshold Ornstein-Uhlenbeck} (TOU) process and
it can be considered as a continuous-time approximation of the
discrete time model \eqref{2-1}. The properties of the MLE and BE of the
threshold $\vartheta $ can be studied with the help of the technique
developed in \cite{Kut04}.

\subsection*{Acknowledgements} We would like to thank J. Dedeker for many helpful
discussions and bringing to our attention the inequality \eqref{DD} and
W. Liu for computational assistance.

\vspace{1.0in}
\hspace{-0.32in}
\begin{tabular}{lll}
Ngai Hang Chan & &Yury A. Kutoyants\\
Department of Statistics & &Laboratoire de Statistique et Processus\\
Chinese University of Hong Kong & &Universit\'e de Maine\\
Shatin, NT & & 72085 Le Mans,  Cedex 9\\
Hong Kong & & France \\
{\sc e-mail}: nhchan@sta.cuhk.edu.hk & &
{\sc e-mail}: kutoyants@univ-lemans.fr\\
\end{tabular}

\begin{center}
\begin{figure}
\includegraphics[height=15cm, width=15cm]{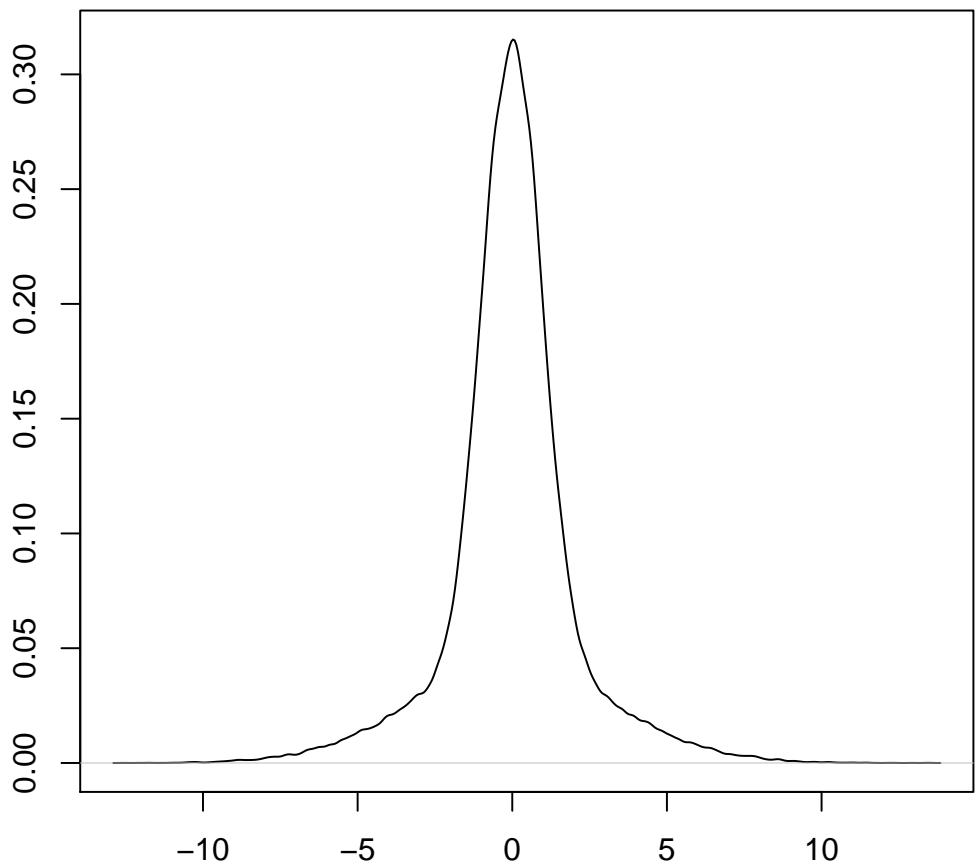}
\caption{The kernel density function of $f(\theta,2)$.}
\end{figure}
\end{center}

\begin{center}
\begin{figure}
\includegraphics[height=15cm, width=15cm]{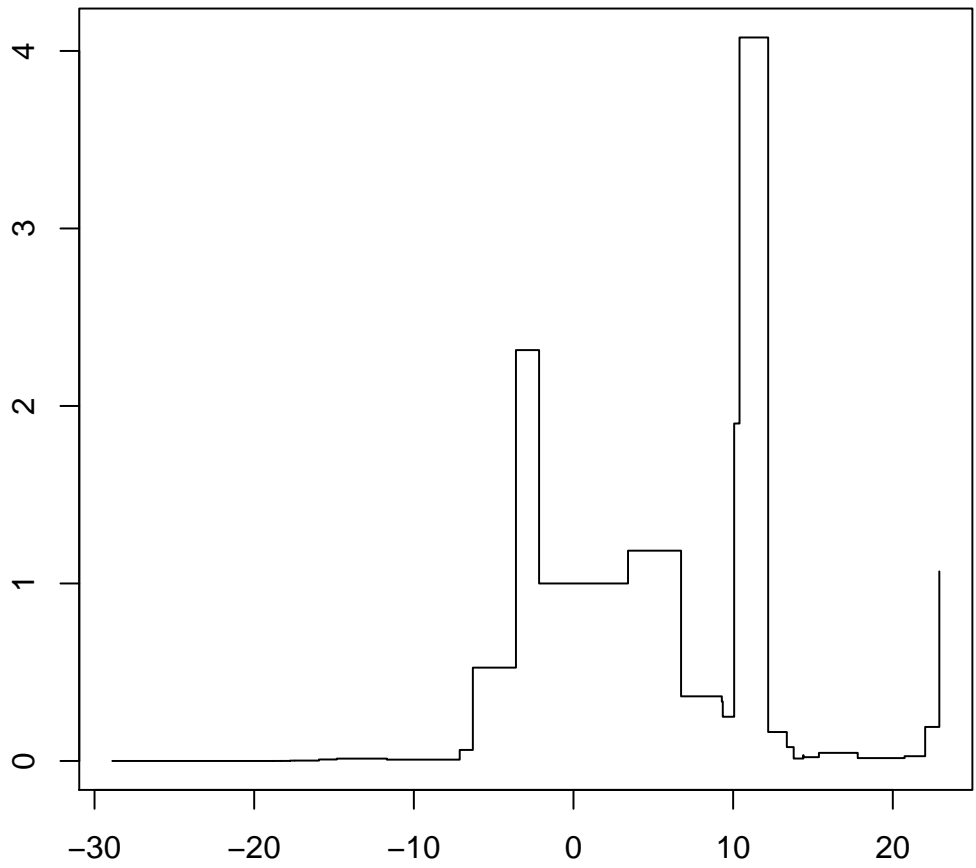}
\caption{A sample path of $Z(u)$.}
\end{figure}
\end{center}

\begin{center}
\begin{figure}
\includegraphics[height=15cm, width=15cm]{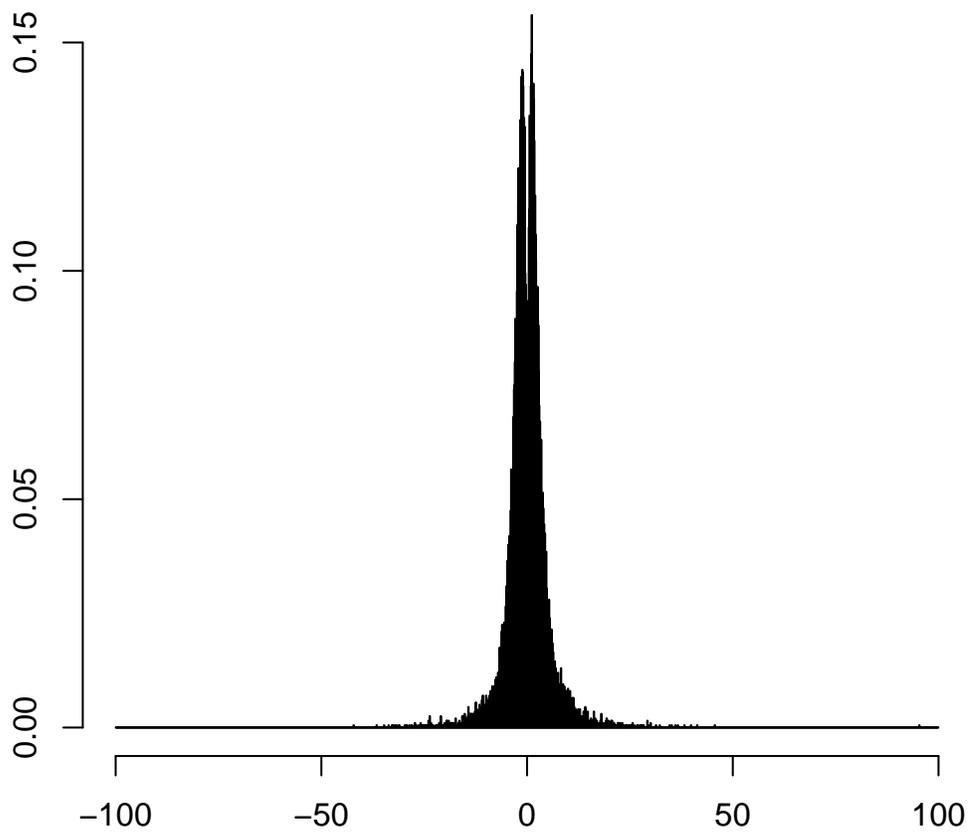}
\caption{Histogram of the MLE $\hat{u}$.}
\end{figure}
\end{center}

\begin{center}
\begin{figure}
\includegraphics[height=15cm, width=15cm]{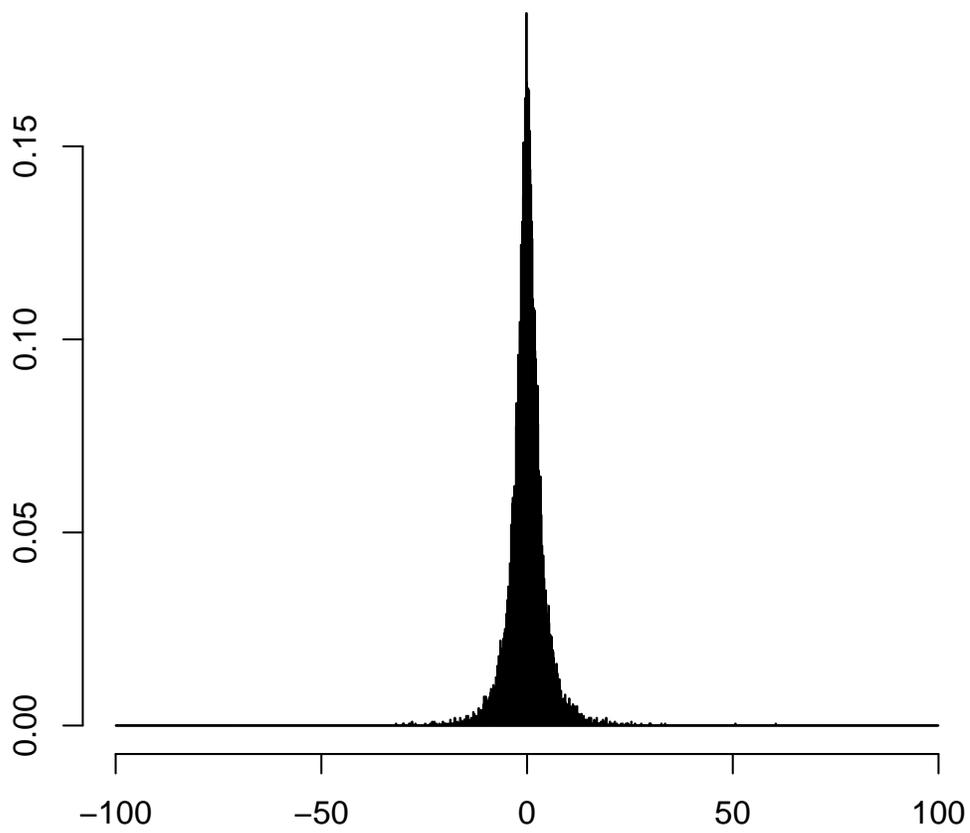}
\caption{Histogram of the Bayesian estimator $\tilde{u}$.}
\end{figure}
\end{center}
\end{document}